\documentclass[a4paper,12pt]{iopart}

\usepackage{topcapt}
\usepackage{amssymb}
\usepackage{amsthm}
\usepackage{amsfonts}
\usepackage{url}
\usepackage{rotating}
\usepackage{adjustbox}
\usepackage{longtable}
\usepackage{tikz}
\usepackage{tcolorbox}
\usepackage{algorithm}
\usepackage{algpseudocode}
\usepackage{pdflscape}
\usepackage{microtype}
\usepackage{caption}
\usepackage{hyperref}

\newcommand{\BE}{\begin{equation}}
\newcommand{\EE}{\end{equation}}

\newtheorem{theorem}{Theorem}

\newcommand{\Sig}{\textstyle{\sum}}
\renewcommand{\Sig}{\raisebox{-0.35ex}{\textrm{\large $\Sigma$}}}
\newcommand{\U}{\raisebox{0.15ex}{\tcbox[on line, arc=2pt,colback=black!10!white,boxrule=0.5pt, boxsep=0.25ex,left=0.5ex,right=0.5ex,top=0.05ex,bottom=0.05ex]{\scriptsize $)$}}}
\renewcommand{\L}{\raisebox{0.15ex}{\tcbox[on line, arc=2pt,colback=black!10!white,boxrule=0.5pt, boxsep=0.25ex,left=0.5ex,right=0.5ex,top=0.05ex,bottom=0.05ex]{\scriptsize $($}}}
\newcommand{\E}{\raisebox{0.15ex}{\tcbox[on line, arc=2pt,colback=black!10!white,boxrule=0.5pt, boxsep=0.25ex,left=0.4ex,right=0.4ex,top=0.5ex,bottom=0.45ex]{\scriptsize $\circ$}}}
\newcommand{\B}{\raisebox{0.15ex}{\tcbox[on line, arc=2pt,colback=black!10!white,boxrule=0.5pt, boxsep=0.25ex,left=0.4ex,right=0.4ex,top=0.5ex,bottom=0.45ex]{\scriptsize $\bullet$}}}
\newcommand{\UB}{\raisebox{0.15ex}{\tcbox[on line,arc=2pt,colback=blue!20!white,boxrule=0.5pt, boxsep=0.25ex,left=0.5ex,right=0.5ex,top=0.05ex,bottom=0.05ex]{\scriptsize)}}}
\newcommand{\LB}{\raisebox{0.15ex}{\tcbox[on line, arc=2pt,colback=blue!20!white,boxrule=0.5pt, boxsep=0.25ex,left=0.5ex,right=0.5ex,top=0.05ex,bottom=0.05ex]{\scriptsize $($}}}
\newcommand{\EB}{\raisebox{0.15ex}{\tcbox[on line, arc=2pt,colback=blue!20!white,boxrule=0.5pt, boxsep=0.25ex,left=0.4ex,right=0.4ex,top=0.5ex,bottom=0.45ex]{\scriptsize $\circ$}}}
\renewcommand{\BB}{\raisebox{0.15ex}{\tcbox[on line, arc=2pt,colback=blue!20!white,boxrule=0.5pt, boxsep=0.25ex,left=0.4ex,right=0.4ex,top=0.5ex,bottom=0.45ex]{\scriptsize $\bullet$}}}

\newcommand{\MZ}{{\mathcal M}}
\newcommand{\C}{{\mathcal C}}
\newcommand{\T}{{\mathcal T}}
\newcommand{\PC}{{\mathcal P}}
\newcommand{\ggs}{\widehat{g}}

\theoremstyle{plain}
\renewcommand{\thetable}{\arabic{table}}

\newcommand{\orcid}[1]{\address{ORCID ID: \href{http://orcid.org/#1}{#1}}}

\begin{document}

\title{The gerrymander sequence, or A348456.}

\author{Anthony J Guttmann}
\address{School of Mathematics and Statistics,
The University of Melbourne,
Vic. 3010, Australia}
\ead{guttmann@unimelb.edu.au}
\orcid{0000-0003-2209-7192}

\author{Iwan Jensen}
\address{College of Science and Engineering, Flinders University at Tonsley,
GPO Box 2100, Adelaide, SA 5001, Australia}
\ead{iwan.jensen@gmail.com}
\orcid{0000-0001-6618-8470}

\begin{abstract}
The {\em gerrymander sequence}, $g_L$, given as A348456 in the OEIS,  counts the number of ways to dissect a $2L \times 2L$ chessboard into two polyominoes of equal area.
Recently Kauers, Koutschan and Spahn announced a significant increase in the length of this sequence from 3  to 7 terms.  We give a further extension to 11 terms, but more significantly prove that the coefficients grow as $\lambda^{4L^2},$ where $\lambda \approx 1.7445498, $ and is equal to the corresponding quantity for self-avoiding walks crossing a square (WCAS), or self-avoiding polygons crossing a square (PCAS). 
These are, respectively, OEIS sequences A007764 and A333323. Thus we have established a close connection between these previously separate problems.

We have also related the sub-dominant behaviour to that of WCAS and PCAS, allowing us to conjecture that the coefficients of the gerrymander sequence grow as $g_L\sim \lambda^{4L^2+2dL+e} \cdot (2L)^h,$ where
$d=-4.04354 \pm 0.0001,$ $e \approx 8$ and $h=0.75 \pm 0.01$, with $h$ almost certainly $3/4$ exactly.

We have also generated 26 terms of the related OEIS sequence A068416, which counts the number of ways to partition a $L\times L$ square into two connected components (not necessarily of equal area). We have thus been able to predict the asymptotic behaviour of this sequence with a satisfying degree of precision. Indeed, it behaves exactly as $L$ times the corresponding coefficient of the generalised gerrymander sequence (defined below).

The improved algorithm we give for counting these sequences is a variation of that which we recently developed for extending a number of sequences for SAWs and SAPs crossing a domain of the square or hexagonal lattices. It makes use of a minimal perfect hash function and in-place memory updating of the arrays for the counts of the number of paths.

\end{abstract}




\section{Introduction}
\label{sec:intro}

The gerrymander sequence, $g_L$, counts the number of ways to dissect a $2L \times 2L$ chessboard into two polyominoes, each of area $2L^2.$ Until very recently only the first three terms of this sequence
were known. Following a challenge from Neil Sloane the sequence was extended by Kauers, Koutschan and Spahn from  3 terms to 7 terms \cite{KKS22}. In this paper (among other calculations) we further extend the gerrymander sequence to 11 terms.

Our approach is slightly more general. Divide an $L \times L$ square into two connected regions {\bf but not necessarily of equal area}.
We refer to these as {\em generalised gerrymander configurations}.
Let $g_{L,k}$ be the number of such configurations with one region having area $k$ (the other having area $L^2-k$). Note that any configuration is counted twice, once for the
region of area $k$ and once for the region of area $L^2-k$. Then, clearly $g_{L,k}$ is symmetric, so that $g_{L,k}=g_{L,L^2-k}.$ 

 Set $i=\lfloor L^2/2\rfloor.$ We then define the {\em generalised gerrymander sequence} as $\ggs_L=g_{L,i}$. 
For $L$ even $\ggs_L$ is twice the gerrymander sequence coefficient,  while for $L$ odd,  it is one of the two (equal)  terms on either side of the half area mark. The sequence $g_{L,k}$ is by definition symmetric and by observation unimodal 
so  $\ggs_L$ is the largest term in the sequence\footnote{We cannot prove unimodality, which is not surprising, as it took 150 years to prove unimodality \cite{PP13} for the much simpler subset of generalised gerrymander configurations in which only steps to the north and east are taken. The area under such a path is given by the $q$-binomial coefficients.}.

It should be noted  that in any generalised gerrymander configuration one (or both) regions has to be a self-avoiding polygon (SAP). 
This is the key to our efficient enumeration of gerrymanders. 

\begin{figure}
\begin{center}
\includegraphics{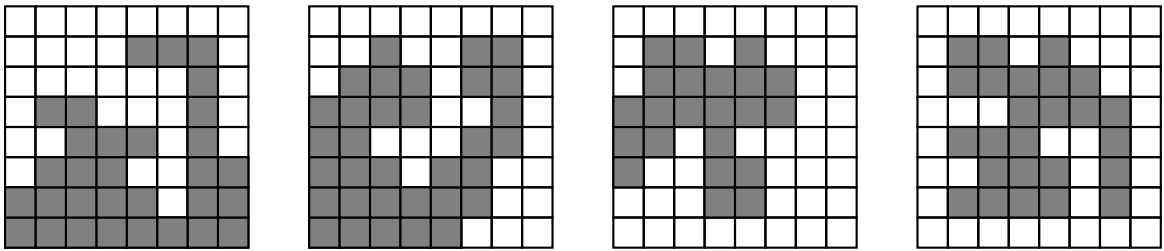}
\end{center}
\caption{\label{fig:gerrySAP} The four cases of self-avoiding polygons (grey regions)  resulting in gerrymander configurations.
}
\end{figure}

In \Fref{fig:gerrySAP} we display examples of the four distinct cases one has to consider, noting that the grey region is a SAP. 
Given the constraint of only two connected regions, it follows that 
the grey region can be chosen so that it contains either zero, one, or two  corners of the square. 

If the grey region contains none of the corners then it is either a SAP not touching any of the sides 
or it is a SAP touching one and only one side of the square. In the first case, shown in panel 4 of  \Fref{fig:gerrySAP}, all cells on the border belongs to the white region. In the second case the part of the SAP along the side must consist of contiguous cells as illustrated in panel 3 of  \Fref{fig:gerrySAP}. All cells on the remaining sides of the square belongs to the white region.

If the grey region contains one corner then the cells along the two sides next to this corner must be  contiguous as shown in panel 2 of \Fref{fig:gerrySAP}. All cells along the two remaining sides (top row and right-most column in our example) belongs to the white region.
Finally when the grey region contains two corners these must lie on the same side of the square (see panel 1 of \Fref{fig:gerrySAP}). All cells on the remaining side (top row in our example) belong to the white region.

By symmetry, configurations shown in panel 1 must be counted twice since the SAP could include the bottom or left side of the square, those in panels 2 and 3 must be counted four times and those in panel 4 must be counted only once. To count generalised gerrymander configurations we therefore need to count SAPs by area (with some constraints which we will detail further in Section~\ref{sec:algo}).

We  define an area generating function (or {\em gerrymander polynomial}) 
\begin{equation}\label{eq:gerrypol}
 G_L(q) = \sum_{k=1}^{L^2-1} g_{L,k} q^k   
\end{equation}
$G_L(q)$ counts all generalised gerrymander configurations 
and is a symmetric polynomial. $G_L(1)/2$ is the number of ways 
to partition the square into two connected regions. We calculated $G_L(q)$ up to $L=22$,
and we extended this up to $L=26$ for the special case of $G_L(1)/2$.
The first 6 terms of $G_L(1)/2$ were obtained by RH Hardin and extended to 14 terms by A Howroyd,
see OEIS sequence A068416.

A {\em $n$-step self-avoiding walk} (SAW) ${\bf \omega}$ on a regular lattice is 
a sequence of {\em distinct} vertices, $\omega_0, \omega_1,\ldots , \omega_n$, 
such that each vertex is a nearest neighbour of its predecessor. SAWs are
considered distinct up to translations of the starting point $\omega_0$.
If  $\omega_0$ and  $\omega_n$ are nearest-neighbours we can form
a closed  $(n+1)$-step self-avoiding polygon (SAP) by adding
an edge between the two end-points.

We recently studied SAWs on an $L \times L$ square lattice, with the walks starting at the south-west corner $(0,0)$ and finishing at the north-east corner $(L,L),$ and constrained within the square (see the first diagram in \Fref{fig:example}), and refer to these as WCAS \cite{GJ22}.
Consider the generating function $C_L(x) = \sum_{n} c_{L,n} x^n,$ where $c_{L,n}$  denotes the number of WCAS of length $n$ on a $L\times L$ square. 

The existence of the limit 
\BE \label{eq:CLlim}
\lim_{L \to \infty} C_{L}(1)^{1/L^2}=\lambda
\EE
was proved in  both \cite{AH78} and \cite{GW90} by different methods, and has recently been given a third proof \cite{SW22}. 
In our recent work \cite{GJ22} we estimated  $\lambda = 1.7445498 \pm 0.0000012.$
We also estimated the sub-dominant terms by finding compelling numerical evidence for the asymptotic behaviour 
\BE \label{eq:CLas}
C_L(1) \sim \lambda^{L^2+bL+c}\cdot L^g,
\EE
 where $b=-0.04354 \pm 0.0001,$ $c=0.5624 \pm 0.0005,$ and $g=0.000 \pm 0.005,$ from which we conjectured that $g=0$, exactly.  This conjectured asymptotic form has received support by a very recent result of Whittington \cite{SW22} who has proved that $C_L(1) = \lambda^{L^2+O(L)},$ a significant improvement on the previous result $C_L(1) = \lambda^{L^2+o(L^2)}.$

For SAPs crossing a square (PCAS) our analysis clearly demonstrated that the growth constant
is the same as for WCAS. We prove this result in Section~\ref{sec:proof}. We also conjectured that the subdominant term $\lambda^b$ is the same as for WCAS, estimated $c \approx -1.197,$ and that the corresponding exponent $g=-\frac12.$

In Section \ref{sec:proof} we prove that the growth constant for the generalised gerrymander sequence equals the growth constant for WCAS.
We use the data generated in this study to estimate the sub-dominant behaviour, and find that $\ggs_L \sim \lambda^{L^2+dL+e} \cdot L^h,$ where $d=-4.04354 \pm 0.0001,$ $e \approx 8$ and $h=0.75 \pm 0.01$ for generalised gerrymanders, and is likely $3/4$ exactly. 
Consequently we find that the coefficients of the gerrymander sequence $g_L$ (OEIS A348456) grow as $g_L\sim \lambda^{4L^2+2dL+e} \cdot (2L)^h$.

Finally, we find $G_L(1)/2 \sim \lambda^{L^2+dL+e} \cdot L^g,$ where now $d=-4.04354 \pm 0.0001,$ $e \approx 8$ and $g=1.75 \pm 0.01,$ and is likely $7/4$ exactly. So this is just $L$ times the corresponding coefficient $\ggs_L.$

In Section~\ref{sec:proof} we present proofs that the growth constants for both PCAS and gerrymanders are equal to the growth constant for WCAS.
In Section~\ref{sec:algo} we give a description of the new and very efficient algorithm we used to calculate the series for generalised gerrymanders.
Section~\ref{sec:ana} contains an  analysis of the generalised gerrymander sequence. 
Section~\ref{sec:conc} contains our conclusions and gives a summary of the estimates we have obtained.

\section{Proof of main results for growth constants} \label{sec:proof} 

\begin{figure}[ht!] 
\centerline{\includegraphics[width=0.9\textwidth,angle=0]{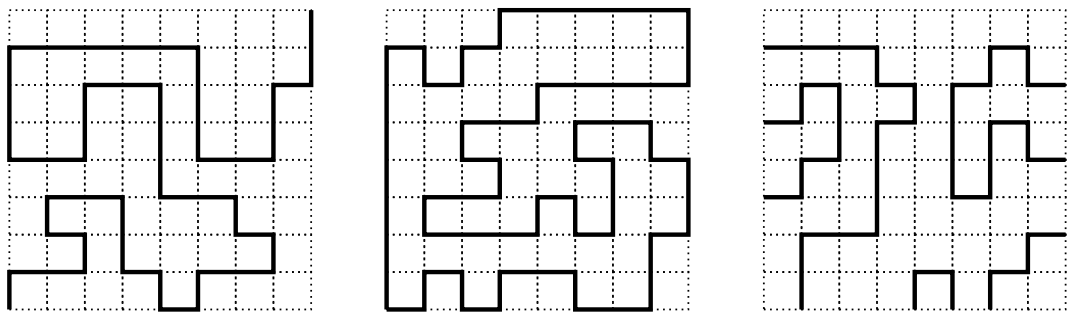} }
\caption{\label{fig:example}
An example of a SAW configuration crossing a square (left panel), traversing
a square from left to right (middle panel) and a cow-patch (right panel).
}
\end{figure}

In  \Fref{fig:example} we show, from left to right, a WCAS, a PCAS, and a cow-patch, so called  as if we colour each domain within the figure alternately black and white, the resulting figure resembles the pattern on the skin of a suitably endowed cow. In \cite{BGJ05} it was proved that WCAS and cow-patches have the same growth constant. We next give a proof that PCAS and WCAS
have the same growth constant.

\begin{figure}
\begin{center}
\includegraphics[width=0.7\textwidth]{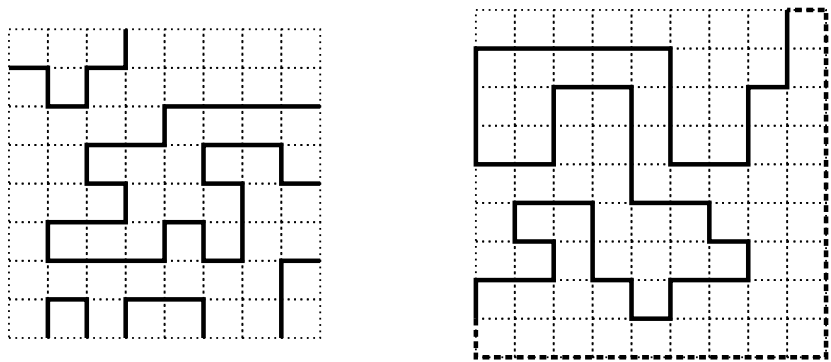}
\end{center}
\caption{\label{fig:bounds}
Illustrations of the constructions used to bound the number of PCAS.
}
\end{figure}

\begin{theorem} \label{th:pcas}
PCAS has the same growth constant as WCAS.
\end{theorem}
\proof 
Consider a PCAS, such as that shown in the middle panel of  \Fref{fig:example}. 
Delete all the bonds along the boundary of the square. 
This transforms the PCAS into a cow-patch pattern (see first panel of \Fref{fig:bounds}). 
Not all cow-patch patterns can be so constructed (the cow-patch in \Fref{fig:example} is an example). 
This shows that the number of PCAS is less than the number of cow-patch patterns in squares of the same size. 

To obtain a bound in the other direction, take a WCAS in an $L \times L$ square, going from $(0,0)$ to $(L,L),$ as shown in the left-most panel of  \Fref{fig:example}. Then add a bond from $(0,0)$ to $ (0,-1),$
 then a sequence of horizontal bonds to $(L+1,-1),$ turn left and add bonds up to $(L+1,L),$ 
then a final bond to $(L,L)$ as shown in the second panel of \Fref{fig:bounds}. 
We have converted a WCAS in a $L \times L$ square to a PCAS in a $(L+1)\times (L+1)$ square. Not
all polygons can be produced this way.
So  $P_{L+1}(1) \ge C_L(1),$ where $P_L(1)$ denotes the number of PCAS in an $L \times L$ square. 

Let $CP_L(1)$ denote the cardinality of the set of cow-patch configurations. We have thus shown that $$C_L(1) \le P_{L+1}(1) \le CP_{L+1}(1).$$ In \cite{BGJ05} it was proved that $\lim_{L \to \infty} C_L(1)^{1/L^2} := \lambda$ exists and is equal to $\lim_{L \to \infty} CP_L(1)^{1/L^2}.$ Hence it follows that $\lim_{L \to \infty} P_L(1)^{1/L^2} = \lambda.$
\qed

Next consider SAPs in an $L \times L$ square that do not necessarily touch any sides
of the square, also known as cycles in the grid graph. Denote the number
of such polygons by ${\mathcal P}_L(1)$. 

\begin{theorem} \label{th:cycle}
${\mathcal P}_L(1)$ grows as WCAS.
\end{theorem}
\proof
Any PCAS is also a cycle so ${\mathcal P}_L(1) \geq P_L(1)$.
If we delete a nominated bond, say the bottom-most, left-most bond, we convert
a cycle into a SAW within a square, but not every SAW can be so constructed.
Let ${\mathcal C}_L(1)$ denote the number of SAWs in a square, then we have

$$ P_L(1) \leq {\mathcal P}_L(1) \leq {\mathcal C}_L(1).$$

We  proved in \cite{GJO22} that ${\mathcal C}_L(1)$ has the same growth constant as
WCAS. As in the previous proof, it then follows that $\lim_{L \to \infty} {\mathcal P}_L(1)^{1/L^2} = \lambda.$
\qed

\begin{theorem} \label{th:GL}
The number of two component partitionings of a square, $G_L(1)/2$, has the same growth constant as WCAS.
\end{theorem}
\proof 
Consider the gerrymander polynomial $G_L(q)$, whose coefficients $g_{L,k}$ are just the number of ways to partition a $L\times L$ square into two connected regions, one of area $k$ and the other of area $L^2-k.$

At least one of the regions has to be a SAP.
So every generalised gerrymander configuration is equivalent to a
SAP in an $L \times L$ square, but not necessarily vice versa. Hence generalised gerrymander configurations are a proper subset of  
SAPs in an $ L \times L$ square.

Any SAP in an $(L-2) \times (L-2)$ square can be surrounded by empty cells on all sides yielding a generalised gerrymander configuration in an $L \times L$ square. 
So generalised gerrymander configurations are a superset of  SAPs on the $(L-2)\times (L-2)$ square.

Let ${\mathcal P}_L(1)$ be the number of 
SAPs in an $L \times L$ square. Then  $P_{L-2}(1) \le G_L(1)/2 \le \mathcal{P}_L(1),$ 
and $G_L(1)/2$ and $C_L(1)$ have the same growth constant. As above, it follows that $\lim_{L \to \infty} (G_L(1)/2)^{1/L^2} = \lambda.$

\qed

\begin{figure}[h]
\begin{center}
\includegraphics[width=0.7\textwidth]{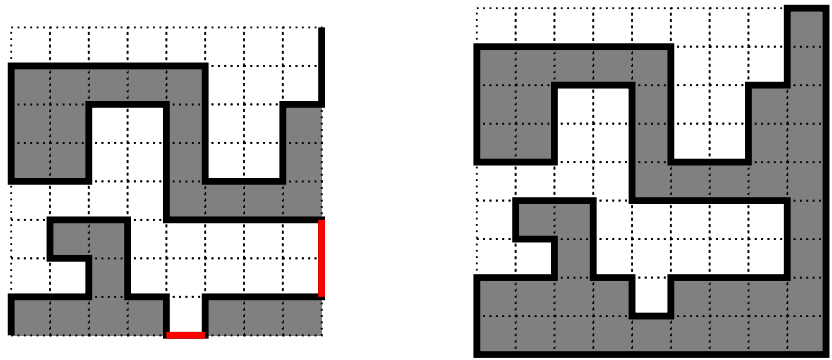}
\end{center}
\caption{\label{fig:wcasa}
The padding of a WCAS with grey cells resulting in a corresponding SAP and hence a connected grey region.
}
\end{figure}

\begin{figure}
\begin{center}
\includegraphics[width=0.7\textwidth]{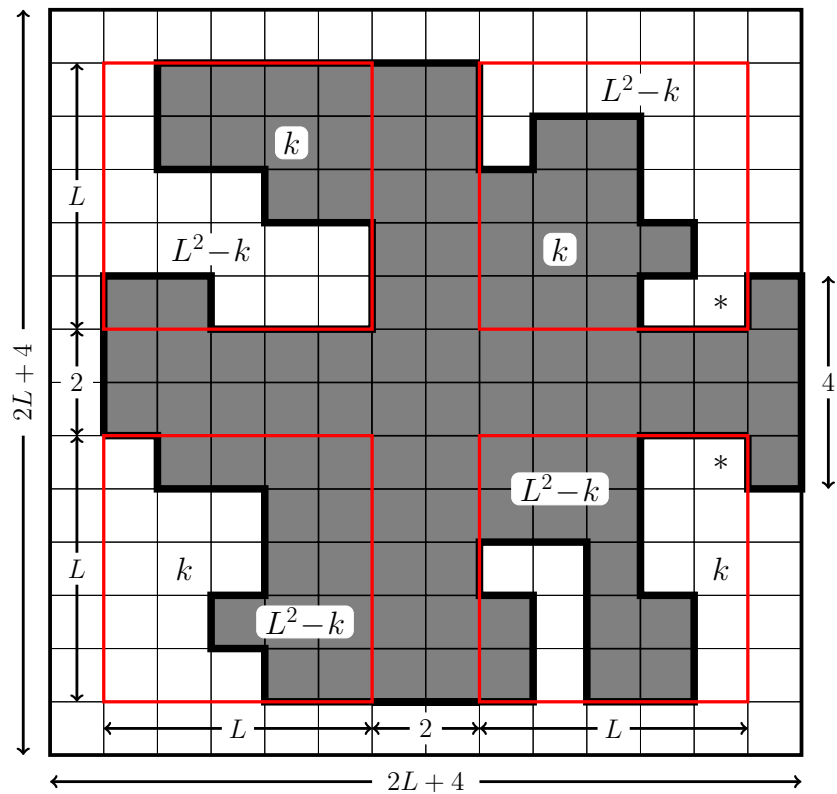}
\end{center}
\caption{\label{fig:lowerbound}
Illustration of the construction used to find a lower bound for the number of gerrymanders in terms of WCAS.
}
\end{figure}

\begin{theorem} \label{th:lower}
The gerrymander sequence has a lower bound,
$g_{L+2} \ge (C_L(1)/L^2)^4$, and $g_L\sim \lambda^{4L^2}$, where $\lambda$
is the growth constant of WCAS.
\end{theorem}
\proof 

Any WCAS from $(0,0)$ to $(L,L)$ divides an $L \times L$ square into two (not necessarily connected) regions, one (grey) of area $n$, the other (white) of area $L^2-n$. 
To see this, navigate the walk from $(0,0)$ to $(L,L),$ shading all squares to the right of the current step -- see the first diagram of \Fref{fig:wcasa}. 
The grey (and indeed white) region will fail to be connected if the walk has steps along the boundary edges from $(0,0)$ to $(L,0),$ or from $(L,0)$ to $(L,L),$ that are not joined to the corners by steps along the boundary (shown in red in  \Fref{fig:wcasa}). 
Similarly, any WCAS touching edges along the other two boundaries  render the white region disconnected. 
The reason being that any walk segment along the boundary renders disjoint the white regions at either end of the segment. 
Now add a row of grey cells immediately below and to the right of the $(L \times L)$ square, plus one additional corner cell where this row and column meet, as shown in the second diagram of \Fref{fig:wcasa}. 
This construction ensures that the grey area is now connected since we can join the end-points of the WCAS by steps along the outside edges of the added cells yielding a SAP as was the case for the construction in \Fref{fig:bounds}. 
Similarly, padding with white cells along the left-
and top-most boundaries will result in a connected white region.

The idea of the proof is to put 4 WCAS together (each within a square of side-length $L,$ with the walk dividing the square into a grey region of area $n$ and a white region of area $L^2-n$) to give a gerrymander in a square of side-length $2L+4$.
Denote the set of all such WCAS as ${\mathcal C}_{L,n}$. 
Then $C_L(1)=\sum_{n=0}^{L^2} |{\mathcal C}_{L,n}|$. 
The maximum (over $n$) of $|{\mathcal C}_{L,n}| \ge C_L(1)/L^2$.
We take this maximum value, whatever it is, say $n=k$, and take any four WCAS in ${\mathcal C}_{L,k}$ and combine them as per \Fref{fig:lowerbound}. 

In the top left corner (inside the square with a red boundary) is a WCAS in ${\mathcal C}_{L,k}$ with area $k$ comprising grey cells. Likewise for the top right corner, though here the walk is reflected so as to cross from $(L,0)$ to $(0,L).$
In the bottom left and bottom right corners are two more WCAS in ${\mathcal C}_{L,k}$, appropriately rotated or reflected, with regions of area $L^2-k$ shaded grey. We then connect these four squares with the grey cross-shaped region of width two. (This adds the required boundary rows and columns of cells to ensure that the grey region is connected). We next add a boundary of white cells of width 1, so that the total width of the square is $2L+4$. This ensures that the white region is connected. 
Finally, we shade grey four squares on the right boundary, as shown. As we show this ensures equal areas of the white and grey regions.

This final shading potentially raises a small pathology, in that the white cell(s) incident on the bottom-right corner can potentially become disjoint. 
This will occur if there is a white cell in the bottom right corner of the WCAS in the first quadrant, marked with an asterisk, (or in the top right corner of the fourth quadrant, similarly marked), and the cell immediately above (below) that corner cell is grey. 
In that case, one must unshade the cell immediately above (below) the corner cell marked with an asterisk, and, correspondingly, shade any white cell which has at least one edge on the walk crossing the square, and which doesn't create a disconnected region. For example, any boundary cell.

The total area of the square is $(2L+4)^2=4L^2+16L+16$. The grey area is $2k+2(L^2-k)+8L+4+4=2L^2+8L+8$,
which is half the total area. So we have constructed a gerrymander.

This construction produces a unique sequence
of gerrymanders. Not all gerrymanders can be so constructed.
So we have, $g_{L+2} \ge C_{L,k}(1)^4 \ge (C_L(1)/L^2)^4$, which gives us our desired lower bound as $C_L(1)^4 \sim \lambda^{4L^2}$. Combined with
the previous result from Theorem \ref{th:GL}, which implies $g_L\leq G_{2L}(1)/2\sim \lambda^{4L^2}$, we
have proved that $\lim_{L \to \infty} g_L^{1/4L^2} =\lambda.$
\qed

\section{Algorithm to enumerate the gerrymander sequence.}
\label{sec:algo}

In Section~\ref{sec:intro} we argued that gerrymander configurations (partitioning the square into two connected regions) can be enumerated by considering the grey SAPs shown in \Fref{fig:gerrySAP}.

The SAPs in panels 1 and 2  can be counted in a single calculation for each $L$ and require counting SAPs in rectangles of size $(L-1)\times L$ with the constraint that the SAP has a single column of cells 
starting from the bottom of the rectangle and with the further constraint that once the SAP leaves the bottom it can never return. 

The SAPs in panel 3 are SAPs in a rectangle of size $(L-2)\times (L-1)$ with the constraint that there is a single column of cells on the left boundary.
Finally, the configurations in panel 4 are unconstrained SAPs in an $(L-2)\times (L-2)$ square. 
These two cases can also be counted in a single calculation as shown later.

To calculate the (generalised) gerrymander sequence we need to count the configurations where each region has equal (or close to) area. To achieve this we actually calculate the complete gerrymander polynomials $G_L(q)$ by enumerating the grey SAPs of \Fref{fig:gerrySAP} by area. 
When combining the counts for the four SAP cases, with appropriate symmetry factors, we get the coefficients $p_{L,k}$, which count the contributions to $G(q)$ from grey SAPs of area $k$ in a $L\times L$ square. 
We then use the coefficients $p_{L,k}$ to calculate $G_L(q)$ since $g_{L,k}=p_{L,k}+p_{L,L^2-k}$, where $p_{L,L^2-k}$ accounts for the contributions from white regions of area $k$. 

\begin{figure}
\centering
\includegraphics[width=0.9\textwidth]{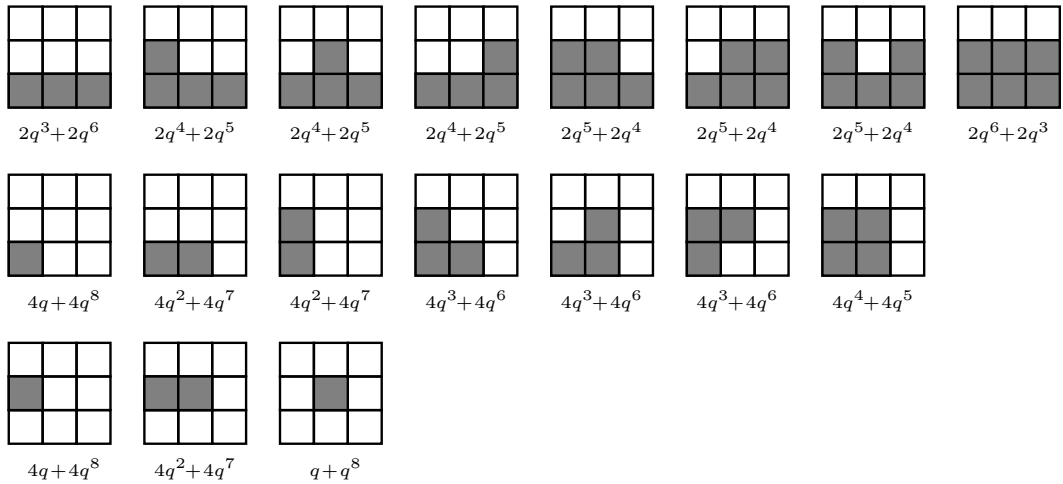}
\caption{\label{fig:gp3} Configurations contributing to $G_3(q)$.}
\end{figure}

In \Fref{fig:gp3} we explicitly list all the configurations 
we must consider to calculate $G_3(q)$. The top row shows the contributions from
configurations of panel 1 in \Fref{fig:gerrySAP}, the middle row arises from 
configurations of panel 2 while the last row arises from configurations in panels 3 and 4.
The contributions to $G_3(q)$ are listed below each panel and the first term comes from 
the grey region while the second term is from the white region. Counting only the
contributions from the grey regions (as done by our algorithms) gives rise to the
polynomial, $9 q +12 q^{2}+14 q^{3}+10 q^{4}+6 q^{5}+2 q^{6}$, which has coefficients
$p_{3,k}$. From these we can then calculate the corresponding gerrymander polynomial,
$$G_3(q)=9 q +12 q^{2}+16 q^{3}+16 q^{4}+16 q^{5}+16 q^{6}+12 q^{7}+9 q^{8}.$$

We implemented the very efficient transfer matrix (TM) algorithm of Iwashita et al. \cite{INK13}  for enumerating SAPs on the square lattice making use of a minimal perfect hash function and in-place memory updating of the arrays for the counts of the number of SAPs. We gave a quite detailed description of the algorithm applied to paths on the hexagonal lattice in our recent paper \cite{GJ22} and most of the considerations from that paper apply to the case of SAPs on the square lattice. So here we will just give a brief description of the main points of the algorithm and point our readers to \cite{GJ22} for further details. 

\begin{figure}
\begin{center}
\includegraphics[scale=0.7]{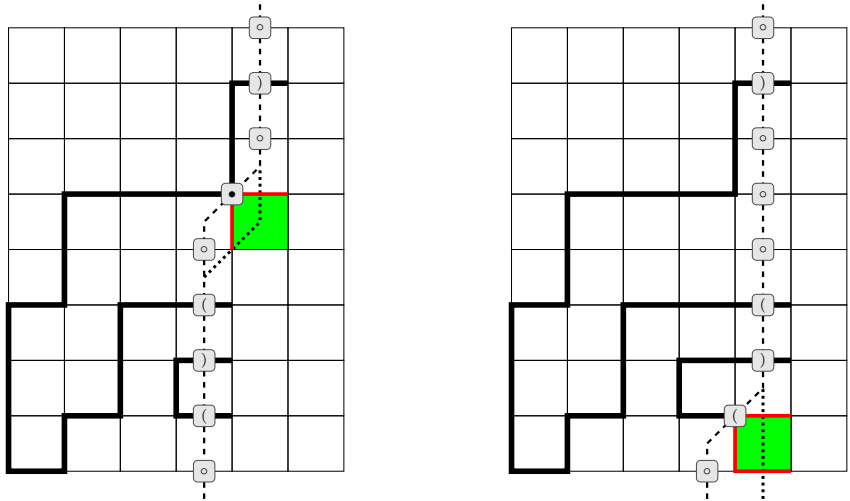}
\end{center}
\caption{\label{fig:TM}
The  basic TM move (left panel) in which the intersection is moved so as to add another cell (coloured green) 
 and two edges (coloured red) to the section of the square already visited. The type of update
 to apply is determined by the states of the edge on the bottom left of the new cell and the vertex
 in the top left corner. The right panel show the final TM move which completes a column of the lattice.}
\end{figure}

If we draw a line across the square  as shown in \Fref{fig:TM} we observe that the partial SAP to the left of the 
intersection consists of arcs connecting two edges on the intersection and a `special' vertex where there is a kink in the intersection. 
On two-dimensional lattices arcs cannot intertwine so each arc end can be assigned a label depending on whether it is the lower or upper end of an arc, and these
labels will form a balanced parenthesis. The vertex at the kink in the intersection can have an additional blocked state if two edges of the partial SAP are incident on the vertex (as shown in the figure). 
We shall refer to the configuration along the intersection as a {\em signature}, denoted by $\Sig$, which can  be represented by a string of states, $\sigma_i$, where

\begin{equation*}
\sigma_i  = \left\{ \begin{array}{rl}
\E &\;\;\; \mbox{empty edge/vertex},  \\ 
\L &\;\;\; \mbox{lower arc end}, \\
\U &\;\;\; \mbox{upper arc end}, \\
\B &\;\;\; \mbox{blocked vertex}. \\
\end{array} \right.
\end{equation*}
\noindent
The partial SAP in \Fref{fig:TM}  has the signature $\Sig = \E\L\U\L\E\B\E\U\E.$ 

\begin{figure}
\begin{center}
\includegraphics[scale=0.7]{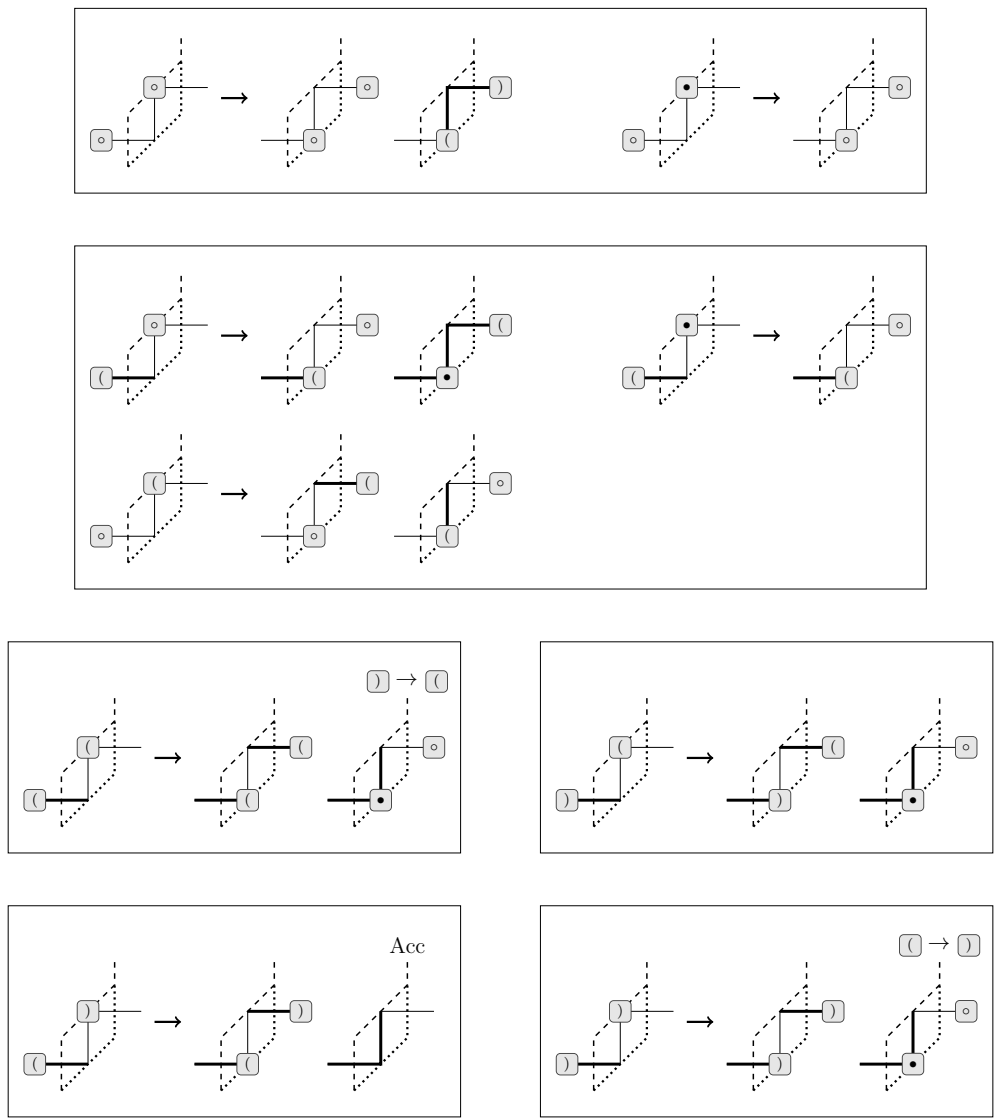}
\end{center}
\caption{\label{fig:upd}
The  possible updates in a TM move with thin edges empty and thick edges occupied by the SAP. The updates are grouped together in such a way as to make 
it possible to use in-place memory updating of the counts.
When two arc ends are joined,  one may have to relabel another arc end as indicated above the update, i.e., for the second transition in the left panel in the third row, two lower arc ends are
joined and the matching upper arc end is relabelled as a lower arc end.}
\end{figure}

For each signature $\Sig$ we simply count the number of partial SAPs, $\C(\Sig,q)$, where $\C(\Sig,q)$ is a polynomial in $q$ such that
the coefficient $p_k$ of $q^k$ equals the number of partial SAPs of area $k$.
SAPs are counted by moving the intersection so as to add a single cell and two edges at a time, as illustrated in \Fref{fig:TM}. For each column of height $W$ the TM move from the left panel of \Fref{fig:TM} is used
$W-1$ times followed by the final move shown in the right panel which completes a column. This sequence of moves is then repeated to construct a rectangle of the required length.

The updating of the counting polynomial $\C(\Sig,q)$ depends on the states of the
bottom edge to the left of the new cell and the topmost vertex of the new cell. The new cell can either belong to a completed SAP or not and
it will belong to a SAP if there is an odd number of occupied edge states below the new cell.

In \Fref{fig:upd} we display the possible local 
`input' states and the `output' states which arise as the kink in the boundary
is propagated by one step (we left out a panel identical to the second panel from the top but involving an upper arc end rather than a lower arc end).  
We shall refer to the signature before the move as the {\em source},  $\Sig_{\rm S}$, and a signature produced as a result of the move as a {\em target},  $\Sig_{\rm T}$. 

It is possible to represent the signatures as  Motzkin paths, which are directed walks  from $(0,0)$ to $(n,0)$  in the first quadrant of the square lattice with step-set 
$\Omega = \{(1,0),(1,1),(1,-1)\}$. The basic mapping from a signature to a Motzkin path is to map \E\ to horizontal steps, \L\ to up steps, and \U\ to down steps. 
Let $\MZ_{n}$ be the set of $n$-step Motzkin paths. Consider enumerations on a 
rectangle of size $W\times L$ ($W\leq L)$.  The set of unblocked signatures is given by the set $\MZ_{W+1}$ since there are $W$ edges and one vertex along the TM intersection. 
A signature can have only a single blocked vertex whose position (in either the source or the target signature) is given by the position of the kink in the TM intersection and the remaining
states form a Motzkin path of length $W$. 
The set of blocked signatures can therefore be represented by the set $\MZ_{W}$.  

Using the mapping of signatures to Motzkin paths one can construct minimal perfect hash functions  $\Phi : \MZ_{W+1} \to \{1,\ldots, |\MZ_{W+1}|\}$
for  unblocked signatures and  $\Psi : \MZ_{W} \to \{1,\ldots, |\MZ_{W}|\}$
for  blocked signatures as described in \cite{INK13,GJ22}.

\begin{algorithm}
\caption{\label{alg:upd}
Update the counts of signatures}
\begin{algorithmic}[1]
 \Procedure {UpdateCounts}{$\Sig_{\rm S}$}
\State $S \gets$ {\sc InputState}($\Sig_{\rm S}$)   \Comment{States of update edges}
\State $\Delta \gets$ {\sc AddArea}($\Sig_{\rm S}$)   \Comment{$\Delta =1$ if new cell in SAP (0 otherwise)}
\If{$S= \EB\EB$}
\State $\Sig_{\rm T} \gets$ {\sc ChangeSignature}($\Sig_{\rm S}$,\LB\UB)   \Comment{Insert new arc}
\State $\Sig_{\rm B} \gets$ {\sc ChangeSignature}($\Sig_{\rm S}$,\BB\EB)   \Comment{Blocked signature}
\State \label{alg:EE-LU}
$\C[\Phi(\Sig_{\rm T}),q] \gets \C[\Phi(\Sig_{\rm T}),q]+\C[\Phi(\Sig_{\rm S}),q]$ \Comment{Update count of $\Sig_{\rm T}$}
\State $\C[\Phi(\Sig_{\rm S}),q] \gets q^\Delta \big(\C[\Phi(\Sig_{\rm S}),q]+\C[\Psi(\Sig_{\rm B}),q]\big)$ \Comment{Update count of $\Sig_{\rm S}$}
\State $\C[\Psi(\Sig_{\rm B}),q] \gets 0$ \Comment{Count of target $\Sig_{\rm B}$ is 0}
\ElsIf{$S= \LB\EB$}
\State $\Sig_{\rm T} \gets$ {\sc ChangeSignature}($\Sig_{\rm S}$,\EB\LB)   
\State $\Sig_{\rm B} \gets$ {\sc ChangeSignature}($\Sig_{\rm S}$,\BB\LB) 
\State $\T(q) \gets C[\Phi(\Sig_{\rm S}),q]$   \Comment{Store count of $\Sig_{\rm S}$}
\State $\C[\Phi(\Sig_{\rm S}),q] \gets q^\Delta\big( \C[\Phi(\Sig_{\rm S}),q]+\C[\Phi(\Sig_{\rm T}),q] +\C[\Psi(\Sig_{\rm B}),q] \big)$ 
\State $\C[\Phi(\Sig_{\rm T}),q] \gets q^{1-\Delta}\C[\Phi(\Sig_{\rm T}),q]$  
\State $\C[\Psi(\Sig_{\rm B}),q] \gets q^{1-\Delta}\T(q)$  \Comment{Count of $\Sig_{\rm B}$ as target}
\ElsIf{$S= \EB\LB$}
\State {\sc Null} \Comment{Do nothing. Processed in previous update}
\ElsIf{$S= \UB\EB$}
\State Same updates as for $\LB\EB$ with $\LB \to \UB$. 
\ElsIf{$S= \EB\UB$}
\State {\sc Null} \Comment{Do nothing. Processed in previous update}
\ElsIf{$S= \LB\LB$}
\State $\Sig_{\rm B} \gets$ {\sc RelabelSignature}($\Sig_{\rm S}$,\EB\BB,\LB)   \Comment{Connect arc ends and relabel}
\State $\C[\Psi(\Sig_{\rm B}),q] \gets \C[\Psi(\Sig_{\rm B}),q]+q^{1-\Delta}\C[\Phi(\Sig_{\rm S}),q]$ 
\State $\C[\Phi(\Sig_{\rm S}),q] \gets q^\Delta\C[\Phi(\Sig_{\rm S}),q]$ 
\ElsIf{$S= \LB\UB$}
\State $\Sig_{\rm T} \gets$ {\sc ChangeSignature}($\Sig_{\rm S}$,\EB\EB)   \Comment{Form closed loop}
 \If{$\Sig_{\rm T} = \E\cdots\E$}          \Comment{Empty signature so valid SAP}
 \State $\PC(q)  \gets \PC(q)  + \C[\Phi(\Sig_{\rm S}),q]$ \Comment{Add to SAP count}
 \EndIf
\State \label{alg:LU} $\C[\Phi(\Sig_{\rm S}),q] \gets q^\Delta\C[\Phi(\Sig_{\rm S}),q]$ 
\ElsIf{$S= \UB\LB$}
\State $\Sig_{\rm B} \gets$ {\sc ChangeSignature}($\Sig_{\rm S}$,\EB\BB)   \Comment{Connect arc ends}
\State $\C[\Psi(\Sig_{\rm B}),q] \gets \C[\Psi(\Sig_{\rm B}),q]+q^{1-\Delta}\C[\Phi(\Sig_{\rm S}),q]$ 
\State $\C[\Phi(\Sig_{\rm S}),q] \gets q^\Delta\C[\Phi(\Sig_{\rm S}),q]$ 
\ElsIf{$S= \UB\UB$}
\State $\Sig_{\rm B} \gets$ {\sc RelabelSignature}($\Sig_{\rm S}$,\EB\BB,\UB)   \Comment{Connect arc ends and relabel}
\State $\C[\Psi(\Sig_{\rm B}),q] \gets \C[\Psi(\Sig_{\rm B}),q]+q^{1-\Delta}\C[\Phi(\Sig_{\rm S}),q]$ 
\State $\C[\Phi(\Sig_{\rm S}),q] \gets q^\Delta\C[\Phi(\Sig_{\rm S}),q]$ 
\EndIf
\EndProcedure
\end{algorithmic}
 \end{algorithm}

\begin{algorithm}
\caption{\label{alg:upd-final}
Update the counts of signatures as final cell added} 
\begin{algorithmic}[1]
 \Procedure {UpdateCounts}{$\Sig_{\rm S}$}
\State $S \gets$ {\sc InputState}($\Sig_{\rm S}$)   \Comment{States of update edges}
\If{$S= \EB\EB$}
\State \label{alg:finEE-LU} $\Sig_{\rm T} \gets$ {\sc ChangeSignature}($\Sig_{\rm S}$,\LB\UB)   \Comment{Insert new arc}
\State $\Sig_{\rm B} \gets$ {\sc ChangeSignature}($\Sig_{\rm S}$,\BB\EB)   \Comment{Blocked signature}
\State $\C[\Phi(\Sig_{\rm T}),q] \gets \C[\Phi(\Sig_{\rm T}),q]+\C[\Phi(\Sig_{\rm S}),q]$ \Comment{Update count of $\Sig_{\rm T}$}
\State $\C[\Phi(\Sig_{\rm S}),q] \gets \C[\Phi(\Sig_{\rm S}),q]+\C[\Psi(\Sig_{\rm B}),q]$ \Comment{Update count of $\Sig_{\rm S}$}
\ElsIf{$S= \LB\EB$}
\State $\Sig_{\rm T} \gets$ {\sc ChangeSignature}($\Sig_{\rm S}$,\EB\LB)   
\State $\Sig_{\rm B} \gets$ {\sc ChangeSignature}($\Sig_{\rm S}$,\BB\LB) 
\State $\C[\Phi(\Sig_{\rm T}),q] \gets  \C[\Phi(\Sig_{\rm T}),q]+\C[\Phi(\Sig_{\rm S}),q]$ 
\State \label{alg:finLE-LE} 
$\C[\Phi(\Sig_{\rm S}),q] \gets q\big(\C[\Phi(\Sig_{\rm T}),q]+\C[\Psi(\Sig_{\rm B}),q]\big)$  
\ElsIf{$S= \EB\LB$}
\State {\sc Null} \Comment{Do nothing. Processed in previous update}
\ElsIf{$S= \LB\LB$}
\State $\Sig_{\rm T} \gets$ {\sc RelabelSignature}($\Sig_{\rm S}$,\EB\EB,\LB)   \Comment{Connect arc ends and relabel}
\State $\C[\Phi(\Sig_{\rm T}),q] \gets \C[\Phi(\Sig_{\rm T}),q]+\C[\Phi(\Sig_{\rm S}),q]$ 
\State $\C[\Phi(\Sig_{\rm S}),q] \gets q\C[\Phi(\Sig_{\rm S}),q]$ 
\ElsIf{$S= \LB\UB$}
\State $\Sig_{\rm T} \gets$ {\sc ChangeSignature}($\Sig_{\rm S}$,\EB\EB)   \Comment{Form closed loop}
 \If{$\Sig_{\rm T} = \E\cdots\E$}          \Comment{Empty signature so valid SAP}
 \State $\PC(q)  \gets \PC(q)  + \C[\Phi(\Sig_{\rm S}),q]$ \Comment{Add to SAP count}
 \EndIf
\State \label{alg:finLU}
$\C[\Phi(\Sig_{\rm S}),q] \gets q\C[\Phi(\Sig_{\rm S}),q]$ 
\EndIf
\EndProcedure
\State After update completed: $\C[\Psi(\Sig_{\rm B}),q] \gets 0, \forall \; \Sig_{\rm B}.$
\end{algorithmic}
 \end{algorithm}

In Algorithms \ref{alg:upd} and \ref{alg:upd-final} we give pseudo code for the updates to the counting polynomials $\C(\Sig,q)$.  
{\sc InputState} simply extracts the states of the edge and vertex involved in the update. 
 {\sc AddArea} determines if the cell added in the TM move adds to the area of the SAP.
{\sc ChangeSignature} changes the states of the input states to those indicated by the two blue tiles. {\sc RelabelSignature} changes the input states to empty states and finds and relabels the matching arc end in those updates where two arc ends are connected in a TM update. 

A technical point should be noted. In the update at line \ref{alg:EE-LU} 
of Algorithm \ref{alg:upd} there should seemingly be a factor $q^{1-\Delta}$ multiplying the term on the right. After all inserting
a new arc changes whether or not the added cell lies within the SAP. The factor is missing at this point because the target signature is processed later as a source and any missing unit of area is added at line \ref{alg:LU}. 
One could rearrange the code and include the processing happening for the case $\LB\UB$ 
when processing the case $\EB\EB$ and then do nothing for the case  $\LB\UB$.

The four cases of SAPs we enumerate can in fact be done in just two separate calculations. 
This requires some minor changes to the code of Algorithms \ref{alg:upd} and \ref{alg:upd-final}
and different ways of initializing $\C(\Sig,q)$ to account for the permitted columns on the 
left-most boundary.

The SAPs in panels 1 and 2 of \Fref{fig:gerrySAP} can be counted in one calculation
on a rectangle of size $(L-1)\times L$. Initially, $\C(\Sig,q)=2q^k$ for signatures 
with a $\L$ at position 0, a $\U$ at position $k$, and all other states $\E$. 
This gives a column of cells  on the left boundary starting from the bottom. 
Next we add $L-1$ columns to the lattice by the TM updates with the addition to 
the SAP count done as: $\PC(q)  \gets \PC(q)  + 2\C[\Phi(\Sig_{\rm S}),q]$.
Then one adds one extra column, but without adding to the SAP count at case $\LB\UB$.
This ensures that all configurations of panel 2 are counted four times. After the
addition of the final column the SAP counts $\PC(q)$ are updated by adding the counts for
the signatures used to initialize $\C(\Sig,q)$, thus counting the configurations in
panel 1 twice. 

One further change to the algorithm is required. Once a SAP has stepped away
from the bottom of the rectangle (there are some empty cells at the bottom of a column) 
it is not allowed to return to the bottom of the rectangle since this would produce 
configurations with more than two connected components. 
This constraint can be easily implemented by some minor changes to Algorithm \ref{alg:upd-final}.
The insertion of a new arc at line \ref{alg:finEE-LU} is not permitted so this line is simply 
removed from the code and the updating at line \ref{alg:finLE-LE} is changed to, 
$\C[\Phi(\Sig_{\rm S}),q] \gets q\big(\C[\Phi(\Sig_{\rm S}),q]+\C[\Psi(\Sig_{\rm B}),q]\big)$,
which prevents the situation where a lower arc end is extended from above to an empty edge on the bottom of the rectangle.

The SAPs in panels 3 and 4 of \Fref{fig:gerrySAP} can be counted in one calculation
on a rectangle of size $(L-2)\times (L-1)$. Initially, $\C(\Sig,q)=4q^{k-j}$ for signatures 
with a $\L$ at position $j$, a $\U$ at position $k>j$, and all other states $\E$. 
This gives a column of cells on the left boundary starting at position $j$ and ending at position $k$.
Next we add $L-2$ columns by TM updates. This counts configurations in panel 3 four times.
The SAPs of panel 4 are just cycles on a $(L-2)\times (L-2)$ square. The counts for
these SAPs are added by putting a few extra lines of code after line \ref{alg:LU} of
Algorithm \ref{alg:upd} and line \ref{alg:finLU} of  Algorithm \ref{alg:upd-final}:

\begin{center}
\begin{algorithmic}
\If{$\Sig_{\rm T} = \E\cdots\E$}          
 \State $\C[\Phi(\Sig_{\rm S}),q]  \gets q + \C[\Phi(\Sig_{\rm S}),q]$ 
 \EndIf
\end{algorithmic}
\end{center}

\noindent
This inserts an arc into the empty signature and it corresponds to starting a
new SAP in an otherwise empty lattice and counting it once.

We used our algorithms to calculate $G_L(q)$ up to $L=22$.

Calculating the sequence $G_L(1)/2$ is a simplification of the algorithms since
we no longer need to keep track of the area of the SAP. That is, the variable $q$ is 
removed and the count $\C(\Sig)$ is just a number and not a polynomial. This makes the
memory use and running time much smaller and we could therefore quite readily extend 
the counts for $G_L(1)/2$ up to $L=26$. 

The coefficients become very large and we deal with this by actually performing all 
calculation modulo various prime numbers and then we reconstruct the actual 
coefficients using the Chinese remainder theorem. For $G_L(q)$ we
used primes of the form $p_i=2^{30}-r_i$ while for $G_L(1)/2$ we used primes
of the form  $p_i=2^{62}-r_i$, with $r_i$ chosen so that we used the largest primes possible.

The algorithms can readily be made parallel using OpenMP for shared memory systems as shown in \cite{GJ22}. 
The most demanding calculation is $L=22$ for SAPs in panel 1 of \Fref{fig:gerrySAP}.
The total memory use was just short of 1TB and required 12 primes. The calculations 
were performed on a system with 2TB of memory and a 48 core 2.6GHz Intel Xeon processor (Icelake).
Using all 48 cores each calculation required 3.26 CPU hours per core and had a wall time of 4.04 hours for a fairly respectable 80\% CPU utilisation rate. 

Similarly for $G_L(1)/2$ the most demanding calculation used about 260GB of memory,
required 8 primes with each using around 1.66 CPU hours per core with a wall time of 1.92 hours
for an 87\% CPU utilisation rate.

The enumeration data for all problems studied in this paper and some of the source code used 
 to calculate the exact coefficients can be found at our GitHub repository \url{https://github.com/IwanJensen/Enumerations/tree/Gerrymander}.

\section{Generated data}

The generalised gerrymander sequence is
{ \scriptsize 
\begin{eqnarray*} 
\fl 0,4, 16, 140, 2804, 161036, 27803749, 14314228378, 21838347160809, 99704315229167288,\\
\fl 1367135978051264146578, 56578717186086829451888706, 7065692298178203128922479762418,\\ 
\fl 2670113158846160742372913777087464324, 3052313665715695874527667027409186333152556,\\
\fl 10576314351887299911761821933016870059157696799590, \\
\fl 111034100174173892447665912670921261073467364516352741228, \\
\fl 3537028455649887297336276306453996860419253673550043079822024000,\\ 
\fl 341733163421465989689352428385746691084586717358593912894419042539233990,\\ 
\fl 100252523974388276666190532080484359784524540996444484455535420554238978388252504,\\ 
\fl 89264965524987466095382312579079040669851719236758669481553261745368627073196518991122982,\\
\fl 241454160053307991366810217012218245577185945222070621347619706812992698342006623035833679925950124.  
\end{eqnarray*}
}
\noindent
This sequence can now be found as entry A358289 in the OEIS \cite{OEIS}. \\
The gerrymander sequence, A348456 is

{\scriptsize
\begin{eqnarray*} 
    \fl 2, 70, 80518, 7157114189, 49852157614583644, 28289358593043414725944353,\\
    \fl 1335056579423080371186456888543732162, 5288157175943649955880910966508435029578848399795,\\ 
    \fl 1768514227824943648668138153226998430209626836775021539911012000,\\ 
    \fl 50126261987194138333095266040242179892262270498222242227767710277119489194126252,\\
    \fl 120727080026653995683405108506109122788592972611035310673809853406496349171003311517916839962975062.
\end{eqnarray*}
}
The sequence $G_L(1)/2,$ previously given in the OEIS \cite{OEIS} to 14 terms as sequence A068416 is\\

\tiny{
\begin{eqnarray*}
  \fl  0, 6, 53, 627, 16213, 1123743, 221984391, 127561384993, 215767063451331, 1082828220389781579, 16209089366362071416785, \\
  \fl  726438398002211876667379681, 97741115155002465272674416929195, 39565596445488219947994403962984729307,\\ 
  \fl  48266553553179571390563558537192580883946581, 177681396812088238354165934687481183466893654956289, \\ 
  \fl 1975937643872352724089992826014929798118573656798037151869, \\
  \fl 66439263265451619293993827233543293728049358568766901433376111533,\\ 
  \fl 6759530908927225810082912389523913516153699624404397503364760209750198965,\\ 
  \fl  2082175573105919327579605927684972257244834479905043487606018480891006067503189435,\\
  \fl 1942921711925290132823971776864238644417661540760987603847662976606820344224527922302926211,\\ 
  \fl 5494467860345971380570753195590132183781118478869271589678543047862011487452888227090933031310123701,\\ 
  \fl 47108220518641518802118971470955381408131440042706598149849757389472954451355452889071492134660868248563001011,\\ 
  \fl 1224939601435250677240944373082090949943585589026155133464810175292336195388826375748636115618176553643756803498606264191,\\  
  \fl 96628963298447280919927953386332031167682264043092747711368137641453332526883992822518172167968007887771992080746478077722442851831,\\ 
  \fl 23130601650038998858250356469453469539571044679175753615930334928835892794769268061796821453552124424327952585726630671538922720993311640741787.
\end{eqnarray*}}
\normalsize

\section{Analysis of sequences.}
\label{sec:ana}
In \cite{GJ22} we gave compelling numerical evidence that PCAS behaved as $P_L(1) \sim \lambda^{L^2+bL+c} \cdot L^g,$ where $\lambda=1.7445498 \pm 0.0000012,$ $b=-0.04354 \pm 0.0001,$ $c \approx -1.197$ and $g = -0.500 \pm 0.005.$

Having shown that the dominant term for generalised gerrymanders is $\lambda^{L^2},$ we make the obvious conjecture that the sub-dominant terms for generalised gerrymanders are similar to those of PCAS, and write $\ggs_L \sim \lambda^{L^2+dL+e} \cdot L^h.$ 

Then \BE \label{eqn:PC}
\frac{P_L(1)}{\ggs_L} \sim \lambda^{\alpha L + \beta} \cdot L^\delta,
\EE
where $\alpha = b-d,$ $\beta=c-e,$ and $\delta=g-h.$

Note that we are analysing the series for {\em generalised gerrymanders.} These will grow like $\lambda^{L^2},$ whereas the gerrymander sequence given by A348456 grows like $\lambda^{4L^2}.$ The asymptotics in that case are clearly stated in the abstract and conclusion.

Equation \ref{eqn:PC} can be rewritten in the more  generic form $\frac{P_L(1)}{\ggs_L}=\widetilde{g}_L \sim F\mu^L \cdot L^\delta,$ where $\mu = \lambda^\alpha,$ and $F=\lambda^\beta.$
Such sequences can be analysed by the ratio method, briefly described in \ref{app:ratio}.

\begin{figure}[ht!] 
\centerline{\includegraphics[width=0.9\textwidth]{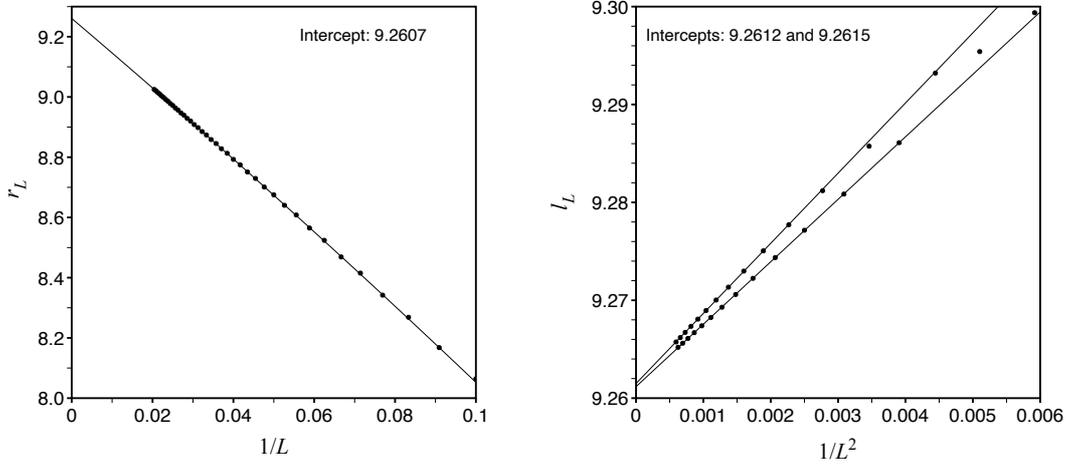}}
 \caption{\label{fig:gerryr} 
Ratios $r_L=\widetilde{g}_L/\widetilde{g}_{L-1}$ plotted against $1/L$ 
and the linear intercepts plotted against $1/L^2$.} 
\end{figure}

We actually extended this sequence by 20 further {\em approximate } terms, using the method of series extension \cite{G16}, briefly described in \ref{app:pred}. The ratios are expected to be accurate to at worst 5 significant digits, which is sufficient accuracy for the ratio method to be used.

In the left panel \Fref{fig:gerryr} we show the ratios, $r_L:=\widetilde{g}_L/\widetilde{g}_{L-1}$, plotted against $1/L.$ In the right panel of \Fref{fig:gerryr} we show the corresponding linear intercepts, 
$l_L=\frac12 \big[L\cdot r_L -(L-2)\cdot r_{L-2}\big]$, which eliminate the term O$(1/L)$ in the ratio sequence, and so should converge faster. We used alternate terms in forming the linear intercepts to minimise a slight parity effect. The solid curves are a quadratic fit using the date from $L=15$ in the case of $r_L$ and a simple linear fit using the data for from $L=20$ in the case of the $l_L$.

\begin{figure}[ht!] 
\centerline{\includegraphics[width=0.9\textwidth]{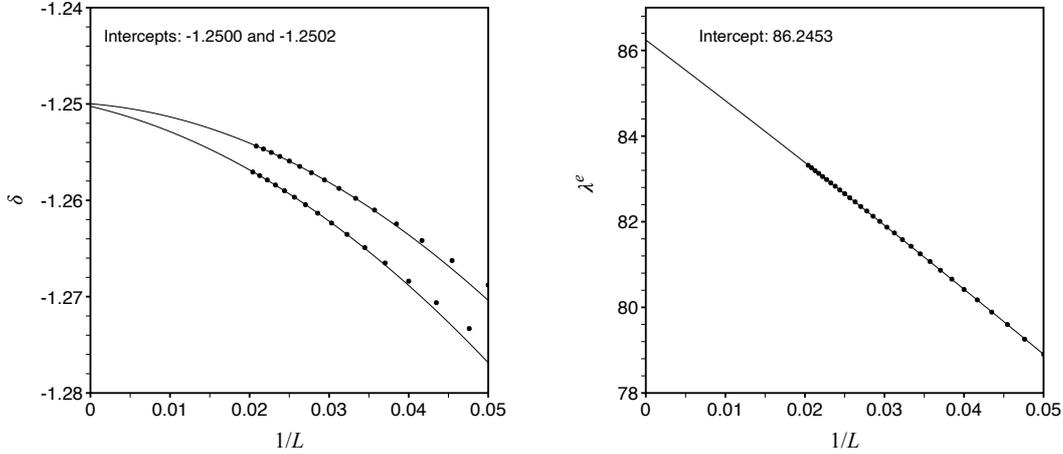}}
 \caption{\label{fig:gerrydelta} 
Estimates of $\delta$ and $\lambda^e$ plotted against $1/L$.} 
\end{figure}

From this figure we estimate the limit as $L \to \infty$ to be around 9.2615. Now this is $\mu=\lambda^\alpha,$ which gives $\alpha \approx 3.9998.$ It is not unreasonable to conjecture that $\alpha=4$ exactly. This then gives $d=-4.04354 \pm 0.0001.$ Assuming $\alpha=4$ exactly, we obtain a refined estimate of $\mu=9.2626123.$ Using this, we estimate the exponent $\delta,$ see eqn. (\ref{A1:3}), and this is shown in the left panel of \Fref{fig:gerrydelta}. 
The solid curves are quadratic fits using the data from $L=30$ (using other starting values gave 
intercepts between $-1.2475$ and $-1.2525$). 
This leads us to conclude $\delta =-1.250 \pm 0.005,$ so that $h=0.750\pm 0.005,$ from which we conjecture $\delta = -5/4$ and $h = 3/4$ exactly. We also eliminated the term O$(1/L)$ which gave an even more convincing plot (not shown), clearly going to a limit very close to $-5/4,$ so we are quite confident in this conjecture.

Finally, we estimate $e$ by extrapolating the sequence $\ggs_L/(\lambda^{L^2+dL} \cdot L^h) \sim \lambda^e,$ using the estimates of $\lambda,$ $d$ and $h$ we have just obtained. This sequence of estimates is shown in the right panel of \Fref{fig:gerrydelta}, and leads us to estimate $\lambda^e \approx 86.2,$ so that $e \approx 8.01.$ Again eliminating the O$(1/n)$ term, we can improve this estimate to $\lambda^e \approx 86.1 \pm 0.2,$ so that $e = 8.01 \pm 0.015.$ Since this estimate is exquisitely sensitive to the estimate of the parameter $d,$ we prefer to quote it as $e \approx 8.$

We next analysed the sequence $G_L(1)/2$ in a similar manner.
Having proved that the dominant term for $G_L(1)$ is $\lambda^{L^2},$ we again conjecture that the sub-dominant terms are similar to those of PCAS, and write $G_L(1)/2 \sim \lambda^{L^2+dL+e} \cdot L^h.$ 

Then \BE \label{eqn:GGF}
\frac{G_L(1)/2}{\lambda^{L^2}} \sim \lambda^{d L + e} \cdot L^h.
\EE

Equation \ref{eqn:GGF} can also be rewritten in the more  generic form $\frac{G_L(1)/2}{\lambda^{L^2}}:=\widetilde{G}_L \sim F\mu^L \cdot L^\delta,$ where $\mu = \lambda^d,$  $F=\lambda^e,$ and $h=\delta.$
This sequence can be analysed similarly to our immediately preceding analysis.

\begin{figure}[ht!] 

\centerline{\includegraphics[width=0.9\textwidth]{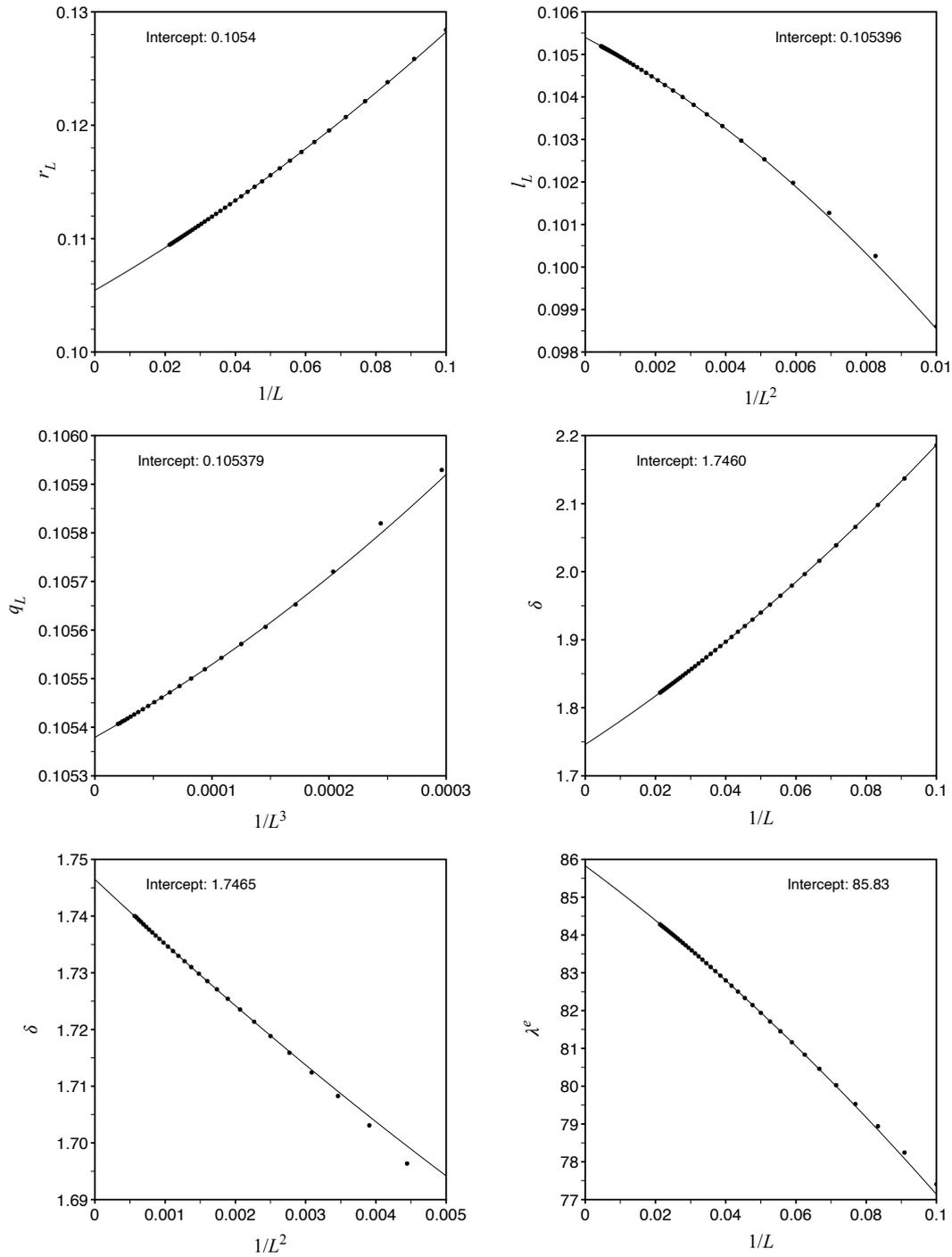}}
 
\caption{\label{fig:gerrypol} 
Ratios $r_L=\widetilde{G}_L/\widetilde{G}_{L-1}$ plotted against $1/L$ (top left panel),
and the corresponding linear (top right panel) and quadratic (middle left panel)
intercepts. Estimates for the exponent $\delta$ (middle right panel) 
and the corresponding liner intercepts (bottom left panel). Estimates of 
the amplitude $\lambda^e$ (bottom right panel).} 

\end{figure}

In the top left panel of \Fref{fig:gerrypol} we show the ratios, $r_L := \widetilde{G}_L/\widetilde{G}_{L-1}\sim \mu(1+\delta/L)$, plotted against $1/L.$ In the top right panel we show the linear intercepts, $l_L := L\cdot r_L - (L-1)\cdot r_{L-1}$, which eliminate the term O$(1/L)$ in the ratio sequence, and so should converge faster.  The middle left panel of  \Fref{fig:gerrypol} shows the quadratic intercepts, $q_L := [L^2\cdot l_L - (L-1)^2\cdot l_{L-1}]/(2L-1)$, which eliminate the term O$(1/L^2).$ These three figures give increasingly precise estimates of $\mu,$ and we conclude that $\mu \approx 0.10538,$ which translates to $d \approx -4.04348.$ Given that the corresponding quantity for the generalised gerrymander sequence that we've just analysed was $-4.04354,$ it seems reasonable to conjecture that this quantity is exactly the same, that is $d = -4.04354,$ so that $\mu \approx 0.1053765.$

Next, to estimate the exponent $\delta,$ we form the sequence $\delta_L := (r_L/\mu -1)\cdot L \sim \delta +O(1/L).$ Estimates of $\delta$ and the linear intercepts, $L\delta_L-(L-1)\delta_{L-1}$, are shown in the middle right and bottom left panels of \Fref{fig:gerrypol}, and lead to the estimate $\delta \approx 1.75=7/4.$

Finally, we estimated the amplitude, just as we did for the previous sequence, by plotting (bottom right panel of \Fref{fig:gerrypol}) the sequence, $G_L(1)/(2\lambda^{L^2+dL} \cdot L^h) \sim \lambda^e,$ using the estimates of $\lambda,$ $d$ and $h$ we have just obtained. This seems to be approaching a limit which we estimate to be $85.8 \pm 0.4,$ so that $e =8.000 \pm 0.008.$ Since this estimate is exquisitely sensitive to the estimate of the parameter $d,$ we prefer to quote it as $e \approx 8,$ exactly the same value as we found for the generalised gerrymander sequence.

\begin{figure}[ht!] 
\centerline{\includegraphics[width=0.9\textwidth]{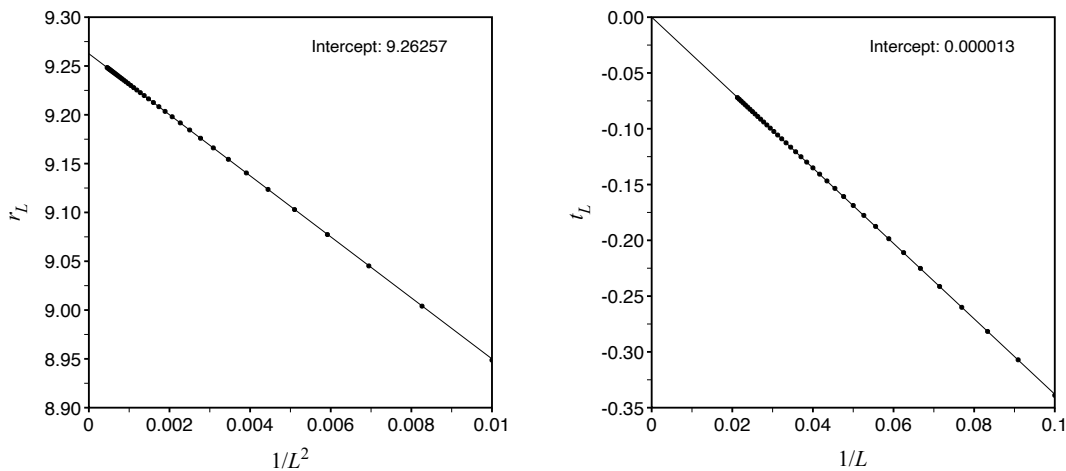}}
\caption{\label{fig:gpolwss}
Ratios, $r_L= \widehat{Sg}_L/\widehat{Sg}_{L-1}$, plotted against $1/L^2$ (left panel) and, $t_L=(r_L/\lambda^4-1)\cdot L$, plotted against $1/L$ (right panel).} 
\end{figure}

For SAWs spanning a square, with coefficients $S_L(1),$ we estimated in \cite{GJ22} that the asymptotics were very similar, with the same exponent $7/4,$ but with growth constant $\mu = \lambda^d$ with $d=-0.04354.$ That is to say, we expect the coefficient ratios to be $\widehat{Sg}_L := S_L(1)/\widehat{g}_L \sim C\cdot \lambda^{4L}.$   To investigate this, we show in the left panel of \Fref{fig:gpolwss} the ratio $r_L = \widehat{Sg}_L/\widehat{Sg}_{L-1},$ which should go to a limit of $\lambda^4 \approx 9.262612$ if our conjectures are correct.
It seems from the figure that this is entirely plausible. To investigate this further, we plotted the quantity $t_L:=(r_L/\lambda^4-1)\cdot L \sim 0$ against $1/L$ in the right panel of  \Fref{fig:gpolwss}, and it can be seen that the data does indeed extrapolate persuasively to zero.

\begin{figure}[ht!] 
\centerline{\includegraphics[width=0.6\textwidth]{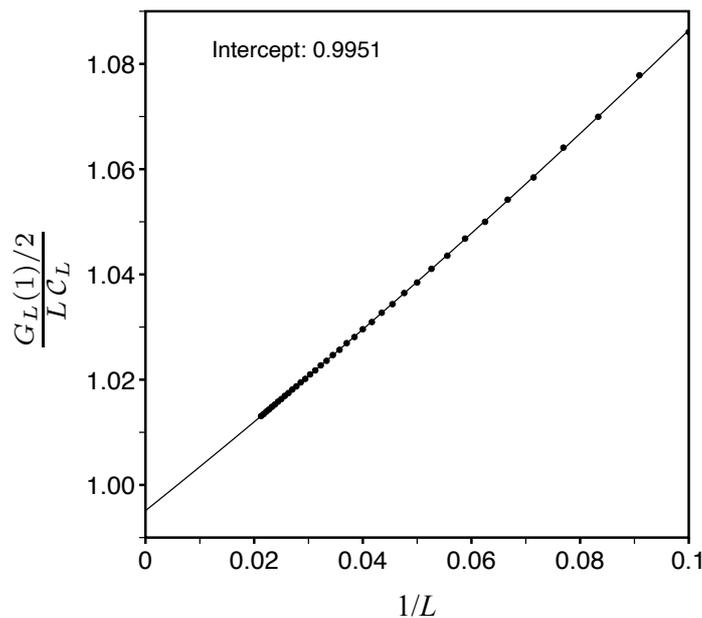}}
\caption{ $\frac{G_L(1)/2}{L\ggs_L}$ plotted against $1/L$.} 
\label{fig:suprat}
\end{figure}

Note that the asymptotics of the two sequences we have analysed are numerically identical, apart from the power of the sub-sub-dominant term, $L^\delta,$ which differs by 1 between the two sequences. To test this observation more precisely, we formed the sequence $$\frac{G_L(1)/2}{ L\cdot \ggs_L} \sim const.$$ 
We show a plot of this ratio in figure \ref{fig:suprat}. This plot is clearly going to a constant value, which appears to be a little below 1. This is abundant support for the conjecture that the asymptotics are identical, apart from a factor of $L.$

\section{Conclusion}
\label{sec:conc}
We have given a new, powerful, algorithm to generate coefficients for the generalised gerrymander sequence $\ggs_L$, and have generated several further terms. 
To be precise, we have generated four further terms in the gerrymander sequence $g_L$, OEIS A348456. We have also generated 26 terms in the sequence $G_L(1)/2$, OEIS A068416, which counts the number of ways to partition a square into two connected regions. 

We have proved that $\ggs_L$ and $G_L(1)/2$ grow with the same dominant behaviour $\lambda^{L^2}$ as do SAWs and SAPs crossing an $L \times L$ square. It then follows that $g_L=\ggs_{2L}/2$ has the dominant behaviour $\lambda^{4L^2}$.

We have used our new data to estimate the sub-dominant terms, and consequently estimate that $\ggs_L \sim \lambda^{L^2+dL+e} \cdot L^h,$ where $d=-4.04354 \pm 0.0001,$ $e \approx 8$ and $h=0.75 \pm 0.01$ for generalised gerrymanders, and the coefficients of the OEIS sequence A348456 grow as $\lambda^{4L^2+2dL+e} \cdot (2L)^h.$

Similarly we find that $G_L(1)/2 \sim \lambda^{L^2+dL+e} \cdot L^h,$ where $d=-4.04354 \pm 0.0001,$ $e \approx 8$ and $h=1.75 \pm 0.01.$ The two sequences $G_L(1)/2$ and $\ggs_L$ behave identically, apart from a factor of $L.$

\section{Acknowledgements}
We would like to thank Nick Beaton for bringing this problem to our attention, Stu Whittington for helpful discussions and a pre-publication version of his paper and Christoph Koutschan and Neil Sloane for careful reading of and useful comments on the manuscript. We particularly wish to thank Andrew Elvey Price for his elucidation of the proof of Theorem \ref{th:lower}.

\newpage
\appendix

\section{Ratio Method}
\label{app:ratio}
 
The ratio method was perhaps the earliest systematic
method of series analysis employed,
and is still the most useful method when only a small number of terms are known. Given a series $\sum c_n z^n,$ which behaves as 
\BE \label{eqn:generic}
F(z)=\sum_n c_nz^n \sim C(1-z/z_c)^{-\gamma},
\EE
it is assumed that $\lim_{n \to \infty} c_n/c_{n-1}$ exists and is equal to the growth constant $\mu = 1/z_c.$ For some combinatorial sequences this has been proved (see e.g.\cite{AM14}), but it is usually just assumed.

From the binomial theorem it follows that $$c_n \sim \frac{C}{\Gamma(\gamma)}\cdot z_c^{-n}
\cdot n^{\gamma-1}.$$From that equation 
it follows that the {\it ratio} of successive terms
\begin{equation} \label{ratios}
r_n = \frac{c_n}{c_{n-1}}=\frac{1}{z_c}\left (1 + \frac{\gamma -1}{n} + {\rm o}\left (\frac{1}{n}\right )\right ).
\end{equation}
 It is then natural to plot the successive ratios $r_n$ against $1/n.$
If the correction terms ${\rm o}(\frac{1}{n})$ can be ignored\footnote{For a purely algebraic singularity, with no confluent terms, the correction term will be ${\rm O}(\frac{1}{n^2}).$}, such a plot will be linear,
with gradient $\frac{\gamma-1}{z_c},$ and intercept $\mu=1/z_c$ at $1/n = 0.$

Linear intercepts $l_n$ eliminate the $O\left ( \frac{1}{n} \right )$ term in eqn. (\ref{ratios}), so in the case of a pure power-law singularity, one has
$$l_n := nr_n - (n-1)r_{n-1} = \mu \left (1+ \frac{c}{n^2} + O\left (\frac{1}{n^3} \right ) \right ).$$ This process can often be iterated, giving quadratic, cubic etc. intercepts.

Various refinements of the method can be readily derived. If the critical point
is known exactly, it follows from eqn. (\ref{ratios}) that estimators of the exponent
$\gamma$ are given by
\BE \label{A1:3}
 \gamma_n: = n(z_c\cdot r_n-1)+1 = \gamma +  {\rm o}(1).
 \EE

If the critical point is not known exactly, one can still estimate the exponent $\gamma.$ From eqn. (\ref{ratios}) it follows that 
\BE \label{eq:exp}
\gamma_n := 1+n^2\left ( 1-\frac{r_n}{r_{n-1}} \right )= \gamma +  {\rm o}(1).
\EE

Similarly, if the exponent $\gamma$ is known, estimators of the growth constant  $\mu$
are given by  $$\mu_n =  \frac{n r_n}{n+\gamma-1}= \mu + {\rm o}(1/n).$$

\section {Differential approximants}
\label{app:da}

The generating
functions  of some problems in enumerative combinatorics are sometimes algebraic, such as that for $Av(1342)$ pattern-avoiding permutations, sometimes D-finite, such as with $Av(12345)$ pattern-avoiding permutations,
sometimes differentially algebraic, and sometimes transcendentally transcendental.
The not infrequent occurrence of D-finite solutions was the origin of the method of {\em differential approximants}, a very successful method of series analysis for analysing power-law singularities \cite{G89}.

The basic idea is to approximate a generating function $F(z)$ by solutions
of differential equations with polynomial coefficients. That is to say, by D-finite ODEs. The singular behaviour
of such ODEs is  well documented
(see e.g. \cite{Forsyth02,Ince27}), and the singular points and
exponents are readily calculated from the ODE. 

The key point for series analysis is that even if {\em globally} the function is not describable by a solution
of such a linear ODE (as is frequently the case) one expects that
{\em locally,} in the
vicinity of the (physical) critical points, the generating
function is still well-approximated by a solution of a linear ODE, when the singularity is a generic power law (\ref{eqn:generic}).

An $M^{th}$-order differential approximant (DA) to a function $F(z)$  is formed by matching
the coefficients in the polynomials $Q_k(z)$ and $P(z)$ of degree $N_k$ and $K$, respectively,
so that the formal solution of the $M^{th}$-order inhomogeneous ordinary differential equation
\BE \label{eq:ana_DA}
\sum_{k=0}^M Q_{k}(z)\left(z\frac{{\rm d}}{{\rm d}z}\right)^k \tilde{F}(z) = P(z)
\EE
agrees with the first $N=K+\sum_k (N_k+1)$ series coefficients of $F(z)$. 

Constructing such ODEs only involves
solving systems of linear equations. The function
$\tilde{F}(z)$ thus agrees with the power series expansion of the (generally unknown)
function $F(z)$ up to the first $N$ series expansion coefficients.
We normalise the DA by setting $Q_M(0)=1,$ thus leaving us with $N$ rather
than $N+1$ unknown coefficients to find. The choice of the differential operator $z\frac{{\rm d}}{{\rm d}z}$ in (\ref{eq:ana_DA}) forces the origin to be a regular singular point. The reason for this choice is that most lattice models with holonomic solutions, for example, the free-energy of the two-dimensional Ising model, possess this property. However this is not an essential choice.

From the theory of ODEs, the singularities of $\tilde{F}(z)$ are approximated by zeros
$z_i, \,\, i=1, \ldots , N_M$ of $Q_M(z),$ and the
associated critical exponents $\gamma_i$ are estimated from the indicial equation. If there is only a single root at $z_i$  this is just
\BE \label{eq:ana_indeq1}
\gamma_i=M-1-\frac{Q_{M-1}(z_i)}{z_iQ_M ' (z_i)}.
\EE
Estimates of the critical amplitude $C$ are rather more difficult to make, involving the integration of the differential approximant. For that reason the simple ratio method approach to estimating critical amplitudes is often used, whenever possible taking into account higher-order asymptotic terms \cite{GJ09}.

Details as to which approximants should be used and how the estimates from many approximants are averaged to give a single estimate are given in \cite{GJ09}. Examples of the application of the method can be found in \cite{G15}. In that work, and in this, we reject so-called {\em defective} approximants, typically those that have a spurious singularity closer to the origin than the radius of convergence as estimated from the bulk of the approximants. Another  method sometimes used is to reject outlying approximants, as judged from a histogram of the location of the critical point (i.e. the radius of convergence) given by the DAs. It is usually the case that such distributions are bell-shaped and rather symmetrical, so rejecting approximants beyond two or three standard deviations is a fairly natural thing to do.

\section{Coefficient prediction}
\label{app:pred}
In analysing combinatorial data, it is often the case that the ratio method and the method of differential approximants  work serendipitously together in many cases, even when one has stretched exponential behaviour, in which case neither method works particularly well in unmodified form. 

To be more precise, the method of differential approximants (DAs)  produces ODEs which, by construction, have solutions whose series expansions agree term by term with the known coefficients used in their construction. Clearly, such ODEs implicitly define {\em all}  coefficients in the generating function, but if $N$ terms are used in the construction of the ODE, all terms of order $z^{N}$ and beyond will be approximate, unless the exact ODE is discovered, in which case the problem is solved, without recourse to approximate methods.

It is useful to construct a number of DAs that use all available coefficients, and then use these to predict subsequent coefficients. Not surprisingly, if this is done for a large number of approximants, it is found that the predicted coefficients of the term of order $z^n,$ where $n > N,$ agree for the first $k(n)$ digits, where $k$ is a decreasing function of $n.$ We take as the predicted coefficients the mean of those produced by the various DAs, with outliers excluded, and as a measure of accuracy we take the number of digits for which the predicted coefficients agree, or the standard deviation. These two measures of uncertainty are usually in reasonable agreement.

Now it makes no logical sense to use the approximate coefficients as input to the method of differential approximants, as we have used the DAs to obtain these coefficients. However there is no logical objection to using the ({\em approximate}) predicted coefficients as input to the ratio method. Indeed, as the ratio method, in its most primitive form, looks at a graphical plot of the ratios, an accuracy of 1 part in $10^4$ or $10^5$ is sufficient, as errors of this magnitude are graphically unobservable.

Ratio methods, and direct fitting methods, by contrast are much more robust. The sort of small error that affects the convergence of DAs would not affect the behaviour of the ratios, or their extrapolants, and would thus be invisible to them. As a consequence, approximate coefficients are just as good as the correct coefficients in such applications, provided they are accurate enough. We re-emphasise that, in the generic situation (\ref{eqn:generic}), ratio type methods will rarely give the level of precision in estimating critical parameters that DAs can give. By contrast, the behaviour of ratios can more clearly reveal features of the asymptotics, such as the fact that a singularity is not of power-law type. This is revealed, for example, by curvature of the ratio plots \cite{G15}.

As an example, consider the OGF for $Av(12453)$ PAPs (see OEIS \cite{OEIS} A116485). This is known to order $x^{38}.$ Let us take the coefficients to order $x^{16}$ and use the method of series extension described above to predict the next 22 ratios, so that we can compare them to the exact ratios. The results, based on 3rd order differential approximants, are shown in \Tref{tab:serpred}. For the first predicted ratio, $r_{18},$ the discrepancy is in the 10th significant digit. For the last predicted ratio, $r_{39}$, the error is in the 5th significant digit. This level of precision is perfectly adequate for ratio analysis.

\begin{table}[htbp!]
   \begin{center}
   \topcaption{Ratios $r_{18}$ to $r_{39}$ actual and predicted from the coefficients of $Av(12453),$ with percentage error shown.}
   \begin{tabular}{|l|l|l|} \hline
   Predicted ratios & Actual ratios & Percentage error\\ \hline
10.654655347& 10.65465504& $4.78 \times 10^{-7}$\\ 
10.828226522& 10.82822539&$1.04 \times 10^{-5}$\\
10.986854456& 10.98685140&$2.79 \times 10^{-5}$\\
11.132386843&11.13238007&$4.78 \times 10^{-5}$\\
11.266382111&11.26636895&$6.08 \times 10^{-5}$\\
11.390163118&11.39013998&$2.03 \times 10^{-4}$\\
11.504857930&11.50482182&$3.14 \times 10^{-4}$\\
11.611441483&11.61138359&$4.99 \times 10^{-4}$\\
11.710743155&11.71066190&$6.94 \times 10^{-4}$\\
11.803496856&11.80338255&$9.68 \times 10^{-4}$\\
11.890333733&11.89017822&$1.31 \times 10^{-3}$\\
12.048402545&12.04814337&$2.15 \times 10^{-3}$\\
12.120553112&12.12022972&$2.67 \times 10^{-3}$\\
12.188650126& 12.18824275&$3.34 \times 10^{-3}$\\
12.252994715&12.25252103&$3.87 \times 10^{-3}$\\
12.313939194&12.31336663&$4.65 \times 10^{-3}$\\
12.371707700&12.37104982&$5.32 \times 10^{-3}$\\
12.426619450&12.42581319&$6.49 \times 10^{-3}$\\
12.478784843&12.47787509&$7.29 \times 10^{-3}$\\
12.528486946&12.52743256&$8.41 \times 10^{-3}$\\
\hline
      \end{tabular}
   \label{tab:serpred}
   \end{center}
\end{table}

In practice we find that the more exact terms we know, the greater is the number of predicted terms, or ratios, that can be predicted.

\end{document}